\documentclass{article}

\newcommand{\be}{\begin{equation}}
\newcommand{\ee}{\end{equation}}
\newcommand{\bea}{\begin{eqnarray}}
\newcommand{\eea}{\end{eqnarray}}
\newcommand{\barray}{\begin{array}}
\newcommand{\earray}{\end{array}}
\newcommand{\pa}{\partial}
\newcommand{\nn}{\nonumber}
\newcommand{\bitem}{\begin{itemize}}
\newcommand{\eitem}{\end{itemize}}
\newtheorem{teo}{Theorem}[section]
\newcommand{\bt}{\begin{teo}}
\newcommand{\et}{\end{teo}}
\newtheorem{Def}{Definition}[section]
\newcommand{\bd}{\begin{Def}}
\newcommand{\ed}{\end{Def}}
\newtheorem{lem}{Lemma}[section]
\newcommand{\bl}{\begin{lem}}
\newcommand{\el}{\end{lem}}
\newtheorem{prop}{Proposition}[section]
\newcommand{\bp}{\begin{prop}}
\newcommand{\ep}{\end{prop}}
\newtheorem{cor}{Corollary}[section]
\newcommand{\bc}{\begin{cor}}
\newcommand{\ec}{\end{cor}}
\newtheorem{ex}{Example}[section]
\newcommand{\bex}{\begin{ex}}
\newcommand{\eex}{\end{ex}}
\newtheorem{rem}{Remark}[section]
\newcommand{\br}{\begin{rem}}
\newcommand{\er}{\end{rem}}

\catcode `\@=11
\@addtoreset{equation}{section}

\begin{document}

\begin{center}
{\Large \textbf{Compatible and almost compatible
\\ pseudo-Riemannian metrics\footnote{This work was
supported by the Alexander von Humboldt Foundation (Germany),
the Russian Foundation for Basic Research (project nos. 99--01--00010 and
 96--15--96027) and
INTAS (project no. 96--0770).}}}
\end{center}

\medskip

\begin{center}
{\large {O. I. Mokhov}}
\end{center}

\bigskip

\section{Introduction. Basic definitions} \label{vved}

In this paper, notions of compatible and almost compatible
Riemannian and pseudo-Riemannian metrics,
which are motivated by the theory of
compatible (local and nonlocal) Poisson structures of hydrodynamic type
and generalize the notion of flat pencil of metrics (this notion
plays an important role in the theory of integrable systems of
hydrodynamic type and the Dubrovin theory of Frobenius manifolds
\cite{[31]}, see also \cite{[30]}--\cite{mokh2}), are introduced
and studied.
Compatible metrics generate compatible Poisson structures of
hydrodynamic type
(these structures are local for flat metrics \cite{[1]}
and they are nonlocal if the metrics are not flat \cite{[38]}--\cite{[39]})
and their investigation is necessary for the theory
of integrable systems of hydrodynamic type.
In ``nonsingular'' case, when eigenvalues of pair of metrics are
distinct, in this paper the complete explicit description of
compatible and almost compatible metrics is obtained.
The ``singular'' case of coinciding eigenvalues of pair of metrics is
considerably more complicated for the complete analysis and
has still not been completely studied.
Nevertheless, the problem on two-component compatible
flat metrics is completely investigated. All such pairs, both
``nonsingular'' and ``singular'', are classified by ours.
In this paper we present the complete description
of nonsingular pairs of two-component flat metrics.
The problems of classification of compatible flat metrics
and compatible metrics of constant Riemannian curvature
are of particular interest, in particular, from the viewpoint
of the theory of Hamiltonian systems of hydrodynamic type.
More detailed classification results for these problems
will be published in another paper.
In the given paper we prove that the approach proposed by
Ferapontov in \cite{[34]} for the study of flat pencils of
 metrics can be also applied (with the corresponding modifications
and corrections) to pencils of metrics of constant Riemannian curvature
and to the general compatible Riemannian and pseudo-Riemannian
metrics. We also correct a mistake which is in \cite{[34]}
in the criterion of compatibility of
local nondegenerate Poisson structures of hydrodynamic type
(or, in other words, compatibility of flat metrics).

We shall use both contravariant metrics $g^{ij} (u)$ with upper
indices,
where $u = (u^1,...,u^N)$ are local coordinates, $1 \leq i, j \leq N$,
and covariant metrics
$g_{ij}(u)$ with lower indices,
$g^{is} (u) g_{sj} (u) = \delta^i_j.$
The indices of coefficients of the Levi--Civita connections
 $\Gamma^i_{jk} (u)$
(the Riemannian connections generated by the corresponding
metrics) and tensors of Riemannian curvature $R^i_{jkl} (u)$
are raised and lowered by the metrics corresponding to them:
$$\Gamma^{ij}_k (u) = g^{is} (u) \Gamma^j_{sk} (u),
 \ \ \ \Gamma^i_{jk} (u) = {1 \over 2} g^{is} (u) \left (
{\pa g_{sk} \over \pa u^j} + {\pa g_{js} \over \pa u^k} -
{\pa g_{jk} \over \pa u^s} \right ),$$
$$R^{ij}_{kl} (u) = g^{is} (u) R^j_{skl} (u), \ \ \
R^i_{jkl} (u) = - {\pa \Gamma^i_{jl} \over \pa u^k}
+ {\pa \Gamma^i_{jk} \over \pa u^l} -
\Gamma^i_{pk} (u) \Gamma^p_{jl} (u)
+ \Gamma^i_{pl} (u) \Gamma^p_{jk} (u).$$

\bd \label{d}
Two contravariant metrics
$g_1^{ij} (u)$ and $g_2^{ij} (u)$ of constant Riemannian curvature
$K_1$ and $K_2$, respectively, are called compatible
if any linear combination of these metrics
\be
g^{ij} (u) = \lambda_1 g_1^{ij} (u) + \lambda_2 g_2^{ij} (u),
\label{comb0}
\ee
where $\lambda_1$ and $\lambda_2$ are arbitrary constants such that
$\det ( g^{ij} (u) ) \not\equiv 0$,
is a metric of constant Riemannian curvature
$\lambda_1 K_1 + \lambda_2 K_2$  and
the coefficients of the corresponding Levi--Civita connections
are related by the same linear formula:
\be
\Gamma^{ij}_k (u) = \lambda_1 \Gamma^{ij}_{1, k} (u) +
\lambda_2 \Gamma^{ij}_{2, k} (u). \label{sv0}
\ee
We shall also say in this case that the metrics
$g_1^{ij} (u)$ and $g_2^{ij} (u)$ form a pencil of
metrics of constant Riemannian curvature.
\ed

\bd \label{d1}
Two Riemannian or pseudo-Riemannian contravariant metrics
$g_1^{ij} (u)$ and $g_2^{ij} (u)$ are called compatible if
for any linear combination of these metrics
\be
g^{ij} (u) = \lambda_1 g_1^{ij} (u) + \lambda_2 g_2^{ij} (u),
\label{comb}
\ee
where $\lambda_1$ and $\lambda_2$ are arbitrary constants
such that $\det ( g^{ij} (u) ) \not\equiv 0$,
the coefficients of the corresponding Levi--Civita connections
and the components of the corresponding tensors of
Riemannian curvature are related by the same linear formula:
\be
\Gamma^{ij}_k (u) = \lambda_1 \Gamma^{ij}_{1, k} (u) +
\lambda_2 \Gamma^{ij}_{2, k} (u), \label{sv}
\ee
\be
R^{ij}_{kl} (u) = \lambda_1 R^{ij}_{1, kl} (u)
+ \lambda_2 R^{ij}_{2, kl} (u).  \label{kr}
\ee
We shall also say in this case that the metrics
$g_1^{ij} (u)$ and $g_2^{ij} (u)$ form a pencil of metrics.
\ed

\bd \label{d2}
Two Riemannian or pseudo-Riemannian contravariant metrics
$g_1^{ij} (u)$ and $g_2^{ij} (u)$ are called almost compatible
if for any linear combination of these metrics (\ref{comb})
relation (\ref{sv}) is fulfilled.
\ed

\bd
Two Riemannian or pseudo-Riemannian metrics
$g_1^{ij} (u)$ and $g_2^{ij} (u)$ are called nonsingular pair
of metrics if the eigenvalues of this pair of metrics, that is,
the roots of the equation
\be
\det ( g_1^{ij} (u) -  \lambda g_2^{ij} (u)) =0,
\ee
are distinct.
\ed

These definitions are motivated by the theory of
compatible Poisson brackets of hydrodynamic type.
In the case if the metrics $g_1^{ij} (u)$ and $g_2^{ij}(u)$ are
flat, that is,
 $R^i_{1, jkl} (u) = R^i_{2, jkl} (u) = 0,$
relation (\ref{kr}) is equivalent to the condition that
an arbitrary linear combination of the flat metrics
$g_1^{ij} (u)$ and $g_2^{ij}(u)$ is also a flat metric
and Definition \ref{d1} is equivalent to the well-known
definition of a flat pencil of metrics or, in other words,
a compatible pair of local nondegenerate
Poisson structures of hydrodynamic type
\cite{[31]} (see also \cite{[30]}--\cite{mokh2}).
In the case if the metrics $g_1^{ij} (u)$ and $g_2^{ij}(u)$ are
metrics of constant Riemannian curvature $K_1$ and $K_2$,
respectively, that is,
$$R^{ij}_{1, kl} (u) = K_1 (\delta^i_k \delta^j_l -
\delta^i_l  \delta^j_k), \ \ \
R^{ij}_{2, kl} (u) = K_2 (\delta^i_k \delta^j_l -
\delta^i_l  \delta^j_k),$$
relation (\ref{kr}) gives the condition that an arbitrary
linear combination of the metrics
 $g_1^{ij} (u)$ and $g_2^{ij}(u)$  (\ref{comb})
is a metric of
constant Riemannian curvature
$\lambda_1 K_1 +
\lambda_2 K_2$ and Definition \ref{d1} is
equivalent to Definition \ref{d} of a pencil of
metrics of constant Riemannian curvature or, in other words,
a compatible pair of the corresponding nonlocal Poisson structures
of hydrodynamic type which were introduced and studied by the author
and Ferapontov in \cite{[38]}.
Compatible metrics of more general type correspond
to compatible pairs of nonlocal Poisson structures
of hydrodynamic type which were introduced and studied by
Ferapontov in \cite{[40]}. They arise, for example,
if we shall use a recursion operator generated by a pair
of compatible Poisson structures of hydrodynamic type and
determining, as is well-known, an infinite sequence of
corresponding Poisson structures.

As was earlier noted by the author in
\cite{[36]}--\cite{mokh2}, condition (\ref{kr}) follows
from condition (\ref{sv}) in the case of certain special
reductions connected with the associativity equations (see also
Theorem \ref{mo1} below).
Of course, it is not by chance.
Under certain very natural and quite general assumptions on
metrics (it is sufficient but not necessary, in particular,
that the eigenvalues of the pair of the metrics under
consideration are distinct),
compatibility of the metrics follows from their almost
compatibility but, generally speaking, in the general case,
it is not true even for flat metrics
 (we shall present here below the corresponding
counterexamples). Correspondingly,
we would like to emphasize that condition (\ref{sv})
which is considerably more simple than condition
(\ref{kr}) ``almost'' guarantees compatibility
of metrics and deserves a separate investigation but, in the general
case, it is necessary to require also
the fulfillment of condition (\ref{kr})
for compatibility of the corresponding Poisson structures of
hydrodynamic type. It is also interesting to find out, does
condition (\ref{kr}) guarantee the fulfillment of condition
(\ref{sv}) or not.

\section{Compatible local Poisson structures of
\\ hydrodynamic type}

The local homogeneous Poisson bracket of the first order, that is,
the Poisson bracket of the form
\be
\{ u^i (x), u^j (y) \} =
g^{ij} (u(x))\, \delta_x (x-y) + b^{ij}_k (u(x)) \, u^k_x \,
\delta (x-y), \label{(1.1)}
\ee
where $u^1,...,u^N$ are local coordinates on a certain smooth
$N$-dimensional manifold $M$, is called a
{\it local Poisson structure of hydrodynamic type} or
{\it Dubrovin--Novikov structure} \cite{[1]}.
Here, $u^i(x),\ 1 \leq i \leq N,$ are functions (fields) of
a single independent variable $x$, and
the coefficients $g^{ij}(u)$ and $b^{ij}_k (u)$ of bracket (\ref{(1.1)})
are smooth functions on $M$.

In other words,
for arbitrary functionals $I[u]$ and $J[u]$
on the space of fields $u^i(x), \ 1 \leq i \leq N,$ a bracket of the form
\be
\{ I,J \} = \int {\delta I \over \delta u^i(x)}
\biggl ( g^{ij}(u(x)) {d \over dx} + b^{ij}_k (u(x))\, u^k_x \biggr )
{\delta J \over \delta u^j(x)} dx
\label{(1.2)}
\ee
is defined and it is required that this bracket is a Poisson bracket,
that is, it is skew-symmetric:
\be
\{ I, J \} = - \{ J, I \}, \label{skew}
\ee
and satisfies the Jocobi identity
\be
\{ \{ I, J \}, K \} + \{ \{ J, K \}, I \} + \{ \{ K, I \}, J \} =0
\label{jacobi}
\ee
for arbitrary functionals $I[u]$, $J[u]$ and $K[u]$.
The skew-symmetry (\ref{skew}) and the Jacobi identity
(\ref{jacobi}) impose very strict conditions on the coefficients
$g^{ij}(u)$ and $b^{ij}_k (u)$ of bracket (\ref{(1.2)}) (these
conditions will be considered below).

The local Poisson structures of hydrodynamic type (\ref{(1.1)})
were introduced and studied by Dubrovin and Novikov
in \cite{[1]}. In this paper, they proposed
a general Hamiltonian approach to the so-called
{\it homogeneous systems of hydrodynamic type}, that is, to
evolutionary quasilinear systems of first-order partial differential
equations
\be
u^i_t = V^i_j (u)\, u^j_x
\label{(1.5)}
\ee
that corresponds to structures (\ref{(1.1)}).

This Hamiltonian approach was motivated by the study
of the equations of Euler hydrodynamics and the
Whitham averaging equations that describe the evolution of
slowly modulated multiphase solutions of partial
differential equations \cite{[3]}.

Local bracket (\ref{(1.2)}) is called
{\it nondegenerate} if
$\det (g^{ij} (u)) \not\equiv 0$.
For general nondegenerate brackets of form (\ref{(1.2)}),
Dubrovin and Novikov proved the following important theorem.

\bt [Dubrovin, Novikov \cite{[1]}] \label{dn}
If $\det (g^{ij} (u)) \not\equiv 0$, then bracket (\ref{(1.2)})
is a Poisson bracket, that is, it is skew-symmetric and satisfies
the Jacobi identity if and only if
\bitem
\item [(1)] $g^{ij} (u)$ is an arbitrary flat
pseudo-Riemannian contravariant metric (a metric of
zero Riemannian curvature),

\item [(2)] $b^{ij}_k (u) = - g^{is} (u) \Gamma ^j_{sk} (u),$ where
$\Gamma^j_{sk} (u)$ is the Riemannian connection generated by
the contravariant metric $g^{ij} (u)$
(the Levi--Civita connection).
\eitem
\et

Consequently, for any local nondegenerate Poisson structure of
hydrodynamic type, there always exist local coordinates
$v^1,...,v^N$ (flat coordinates of the metric $g^{ij}(u)$) in which
the coefficients of the brackets are constant:
\be
\widetilde g^{ij} (v)
= \eta^{ij} = {\rm \ const}, \ \
\widetilde \Gamma^i_{jk} (v) = 0, \ \
\widetilde b^{ij}_k (v) =0,
\ee
that is, the bracket has the constant form
\be
\{ I,J \} = \int {\delta I \over \delta v^i(x)}
 \eta^{ij} {d \over dx}
{\delta J \over \delta v^j(x)} dx,
\label{(1.6)}
\ee
where $(\eta^{ij})$ is a nondegenerate symmetric constant matrix:
\be
\eta^{ij} = \eta^{ji}, \ \ \eta^{ij} = {\rm const},
\  \ \det \, (\eta^{ij}) \neq 0.\nn
\ee

On the other hand, as early as 1978, Magri proposed
a bi-Hamiltonian approach to the integration of nonlinear
systems \cite{[11]}. This approach demonstrated
that the integrability is closely related to the
bi-Hamiltonian property, that is, to the property
of a system to have two compatible Hamiltonian
representations. As was shown by Magri in \cite{[11]},
compatible Poisson brackets generate integrable hierarchies of
systems of differential equations.
Therefore, the description of compatible
Poisson structures is very urgent
and important problem in the theory of
integrable systems.
In particular, for a system, the bi-Hamiltonian property
generates recurrent relations for the conservation laws of
this system.

Beginning from \cite{[11]}, quite extensive literature
(see, for example,
\cite{[12]}--\cite{[16]}
and the necessary references therein) has been
devoted to the bi-Hamiltonian approach and to
the construction of compatible Poisson structures
for many specific important equations of
mathematical physics and field theory.
As far as the problem of description
of sufficiently wide classes
of compatible Poisson structures of defined special types is
concerned, apparently the first such statement was considered
in \cite{[17]}, \cite{[18]} (see also \cite{[19]}, \cite{[20]}).
In those papers, the present author posed and completely solved
the problem of description of all compatible local scalar first-order
and third-order Poisson brackets, that is, all Poisson brackets given
by arbitrary scalar first-order and third-order ordinary
differential
operators.
These brackets generalize the well-known compatible pair of
the Gardner--Zakharov--Faddeev bracket
\cite{[21]}, \cite{[22]}
(first-order bracket) and the Magri bracket \cite{[11]}
(third-order bracket) for the Korteweg--de Vries equation.

In the case of homogeneous systems of hydrodynamic type, many integrable
systems possess compatible Poisson structures of hydrodynamic type.
The problems of description of these structures for particular systems and
numerous examples were considered in many papers
(see, for example, \cite{[23]}--\cite{[29]}).
In particular, in \cite{[23]}
Nutku studied a special class of compatible
two-component Poisson structures of hydrodynamic type and
the related bi-Hamiltonian hydrodynamic systems.
In \cite{[28]} Ferapontov classified all two-component
homogeneous systems of hydrodynamic type possessing
three compatible local Poisson structures of
hydrodynamic type.

In the general form, the problem of description of
flat pencil of metrics (or, in other words,
compatible nondegenerate local Poisson structures of hydrodynamic type)
was considered by Dubrovin in
\cite{[31]}, \cite{[30]}
in connection with the construction of important examples of
such flat pencils of metrics, generated by natural pairs of
flat metrics on the spaces of orbits of Coxeter groups and on other
Frobenius manifolds and associated with the
corresponding quasi-homogeneous solutions of
the associativity equations.
In the theory of Frobenius manifolds introduced and studied by
Dubrovin
\cite{[31]}, \cite{[30]}
(they correspond to two-dimensional topological field theories),
a key role is played by flat pencils of metrics, possessing
a number of special additional (and very strict) properties
(they satisfy the so-called quasi-homogeneity property).
In addition, in \cite{[33]} Dubrovin proved that the theory of
Frobenius manifolds is equivalent to the theory
quasi-homogeneous compatible nondegenerate Poisson structures of
hydrodynamic type. The general problem of compatible
nondegenerate local Poisson structures was also considered by Ferapontov
in \cite{[34]}.

The author's papers
\cite{[35]}--\cite{mokh2}
are devoted to the general problem of classification of
local Poisson structures of hydrodynamic type, to integrable
nonlinear systems which describe such
compatible Poisson structures and to special reductions connected
with the associativity equations.

\bd [Magri \cite{[11]}] \label{dm}
Two Poisson brackets $\{ \ , \ \}_1$
and $\{ \ , \ \}_2$ are called {\it compatible} if
an arbitrary linear combination of these Poisson brackets
\be
\{ \ , \ \} = \lambda_1 \, \{ \ , \ \}_1 +
\lambda_2 \, \{ \ , \ \}_2, \label{magri}
\ee
where $\lambda_1$ and $\lambda_2$ are arbitrary constants,
is also always a Poisson bracket.
In this case, one can say also that the brackets $\{ \ , \ \}_1$ and
$\{ \ , \ \}_2$ form a pencil of Poisson brackets.
\ed

Correspondingly, the problem of description
of compatible nondegenerate local Poisson structures
of hydrodynamic type is pure differential-geometric
problem of description of flat pencils of metrics
(see \cite{[31]}, \cite{[30]}).

In \cite{[31]}, \cite{[30]}
Dubrovin presented all the tensor relations
for the general flat pencils of metrics.
First, we introduce the necessary notation.
Let $\nabla_1$ and $\nabla_2$
be the operators of covariant differentiation given by the Levi--Civita
connections $\Gamma^{ij}_{1,k} (u)$
and $\Gamma^{ij}_{2,k} (u)$, generated by the metrics $g^{ij}_1 (u)$
and $g^{ij}_2 (u)$, respectively. The indices
of the covariant differentials are raised and lowered by
the
corresponding metrics: $\nabla^i_1= g^{is}_1 (u) \nabla_{1,s}$,
$\nabla^i_2=g^{is}_2 (u) \nabla_{2,s}$.
Consider the tensor
\be
\Delta ^{ijk} (u) = g^{is}_1 (u) g^{jp}_2 (u)
\left (\Gamma^k_{2, ps} (u)
- \Gamma^k_{1, ps} (u) \right ),  \label{(2.3)}
\ee
introduced by Dubrovin in \cite{[31]}, \cite{[30]}.

\bt [Dubrovin \cite{[31]}, \cite{[30]}] \label{dub1}
 If metrics
$g^{ij}_1 (u)$ and $g^{ij}_2 (u)$ form a flat pencil,
then there exists a vector field $f^i (u)$ such that the tensor
$\Delta ^{ijk} (u)$ and the metric $g^{ij}_1 (u)$ have the form
\be
\Delta ^{ijk} (u) = \nabla^i_2 \nabla^j_2 f^k (u),
\label{(2.4)}
\ee
\be
g^{ij}_1 (u) = \nabla^i_2 f^j (u) + \nabla^j_2 f^i (u) +
c g^{ij}_2 (u), \label{(2.5)}
\ee
where $c$ is a certain constant,
and the vector field $f^i (u)$ satisfies the equations
\be
\Delta^{ij}_s (u) \Delta^{sk}_l (u) =
\Delta^{ik}_s (u) \Delta^{sj}_l (u), \label{(2.6)}
\ee
where
\be
\Delta^{ij}_k (u) =g_{2,ks} (u) \Delta ^{sij} (u)
= \nabla_{2,k} \nabla^i_2 f^j (u), \label{(2.7)}
\ee
and
\be
(g^{is}_1 (u) g^{jp}_2 (u) - g^{is}_2 (u) g^{jp}_1 (u))
\nabla_{2,s} \nabla_{2,p} f^k (u) =0. \label{(2.8)}
\ee
Conversely, for the flat metric $g^{ij}_2 (u)$ and the vector field
$f^i (u)$ that is a solution of the system of equations
(\ref{(2.6)}) and (\ref{(2.8)}),
the metrics $g^{ij}_2 (u)$ and (\ref{(2.5)}) form a flat pencil.
\et

The proof of this theorem immediately follows from the relations that
are equivalent to the fact that the metrics $g^{ij}_1 (u)$
and $g^{ij}_2 (u)$ form a flat pencil and are considered in
flat coordinates of the metric $g^{ij}_2 (u)$ \cite{[31]}, \cite{[30]}.

In my paper \cite{[35]}, an explicit and simple {\em criterion}
of compatibility for two Poisson structures of
hydrodynamic type is formulated, that is, it is shown
what explicit form is sufficient and necessary for the
Poisson structures of hydrodynamic type to be
compatible.

For the moment, we are able to formulate such explicit general criterion
only namely in terms of Poisson structures but not in terms of
metrics as in Theorem \ref{dub1}.

\bl [\cite{[35]}]  \label{lem1}
{\bf (An explicit criterion of
compatibility for Poisson structures of
hydrodynamic type)}
Any local Poisson structure of hydrodynamic type
$\{ I, J \}_2$  is compatible with the constant nondegenerate
Poisson bracket (\ref{(1.6)}) if and only if it has the form
\be
\{ I, J \}_2 =
\int {\delta I \over \delta v^i (x)} \biggl (
\biggl ( \eta^{is} {\partial h^j \over \partial v^s}
+ \eta^{js} {\partial h^i \over \partial v^s} \biggr )
{d \over dx} + \eta^{is}
{\partial^2 h^j \over \partial v^s \partial v^k}
v^k_x \biggr ) {\delta J \over \delta v^j (x)} dx, \label{(2.9)}
\ee
where $h^i (v), \ 1 \leq i \leq N,$ are smooth
functions defined on a certain neighbourhood.
\el

We do not require in Lemma \ref{lem1} that
the Poisson structure of hydrodynamic type
$\{ I, J \}_2$ is nondegenerate.
Besides, it is important to note that this statement is local.

In 1995, in the paper \cite{[34]}, Ferapontov
proposed an approach to the problem on flat pencils of metrics,
which is motivated by the theory of recursion operators,
and formulated the following theorem as a criterion of
compatibility of nondegenerate local Poisson structures
of hydrodynamic type:

\bt [\cite{[34]}] \label{tfer}
Two local nondegenerate Poisson structures of hydrodynamic type
given by flat metrics $g_1^{ij}(u)$ and $g_2^{ij}(u)$
are compatible if and only if the Nijenhuis tensor of the affinor
$v^i_j (u) = g_1^{is} (u) g_{2, sj} (u)$ vanishes,
that is,
\be
N^k_{ij} (u) = v^s_i (u) {\pa v^k_j \over \pa u^s}
- v^s_j (u) {\pa v^k_i \over \pa u^s}  +
v^k_s (u) {\pa v^s_i \over \pa u^j} -
v^k_s {\pa v^s_j \over \pa u^i} =0.
\ee
\et

Besides, it is noted in the remark in \cite{[34]} that
if the spectrum of $v^i_j (u)$ is simple,
the vanishing of the Nijenhuis tensor implies
the existence of coordinates $R^1,...,R^N$ for which
all the objects $v^i_j (u)$, $g_1^{ij} (u)$, $g_2^{ij} (u)$
become diagonal.
Moreover, in these coordinates the
$i$th eigenvalue of $v^i_j (u)$ depends only on the
coordinate $R^i$.
In the case when all the eigenvalues are nonconstant,
they can be introduced as new coordinates.
In these new coordinates
$\tilde v^i_j (R) = {\rm diag \ } (R^1,...,R^N)$,
$\tilde g_2^{ij} (R) = {\rm diag \ } (g^1 (R),...,g^N (R))$,
$\tilde g_1^{ij} (R) = {\rm diag \ } (R^1 g^1 (R),..., R^N g^N (R))$.

Unfortunately, in the general case the Theorem \ref{tfer} is not true
and, correspondingly, it is not a criterion of
compatibility of flat metrics.
Generally speaking, compatibility of flat metrics does not follow
from the vanishing of the corresponding Nijenhuis tensor.
In this paper we shall present the corresponding counterexamples.
In the general case, as it will be shown here,
the Theorem \ref{tfer} is actually a criterion of
almost compatibility of flat metrics that
does not guarantee compatibility of the corresponding
nondegenerate local Poisson structures of hydrodynamic type.
But if the spectrum of $v^i_j (u)$ is simple, that is, all the
eigenvalues are distinct, then we proves here that
the Theorem \ref{tfer} is not only true but also can be
essentially  generalized for the case of arbitrary
compatible Riemannian or pseudo-Riemannian metrics,
in particular, for the especially important cases in
the theory of systems of hydrodynamiic type, namely,
the cases of metrics of constant Riemannian curvature
or the metrics generating the general nonlocal Poisson
structures of hydrodynamic type.
Namely, the following theorem is true.
\bt \label{tmo1}
\bitem
\item[1)] If for any linear combination (\ref{comb})
of two metrics $g_1^{ij} (u)$ and $g_2^{ij} (u)$ condition (\ref{sv})
is fulfilled, then the Nijenhuis tensor of the affinor
$$v^i_j (u) = g_1^{is} (u) g_{2, sj} (u)$$
vanishes.
Thus, for compatible and almost compatible metrics,
the corresponding Nijenhuis tensor always vanishes.
\item[2)] If a pair of metrics
$g_1^{ij} (u)$ and $g_2^{ij} (u)$ is nonsingular, that is,
the roots of the equation
\be
\det ( g_1^{ij} (u) -  \lambda g_2^{ij} (u)) =0,
\ee
are distinct, then it follows from the vanishing
of the Nijenhuis tensor of the affinor
$v^i_j (u) = g_1^{is} (u) g_{2, sj} (u)$ that the metrics
$g_1^{ij} (u)$ and $g_2^{ij} (u)$ are compatible.
Thus, a nonsingular pair of metrics is compatible
if and only if the metrics are almost compatible.
\eitem
\et

\section{Almost compatible metrics and Nijenhuis tensor} \label{sect1}

Let us consider two contravariant Riemannian or pseudo-Riemannian
metrics
$g_1^{ij} (u)$  and  $g_2^{ij} (u)$, and also
the corresponding coefficients of the Levi--Civita connections
$\Gamma^{ij}_{1, k} (u)$ and $\Gamma^{ij}_{2, k} (u)$.

We introduce the tensor

\be
M^{ijk} (u) = g_1^{is} (u) \Gamma_{2, s}^{jk} (u) -
g_2^{js} (u) \Gamma_{1, s}^{ik} (u) -
g_1^{js} (u) \Gamma_{2, s}^{ik} (u) +
g_2^{is} (u) \Gamma_{1, s}^{jk} (u).
\ee

It follows from the following representation
that $M^{ijk} (u)$ is actually a tensor:
\bea
&&
M^{ijk} (u) = g_1^{is} (u) g_2^{jp} (u) ( \Gamma_{2, ps}^k (u) -
 \Gamma_{1, ps}^k (u) ) - \nn\\
&&
g_1^{js} (u) g_2^{ip} (u) (\Gamma_{2, ps}^k (u) -
 \Gamma_{1, ps}^k (u)) =
\Delta^{ijk} (u) - \Delta^{jik} (u).
\eea

\bl \label{l1}
The tensor $M^{ijk} (u)$ vanishes if and only if
the metrics $g_1^{ij} (u)$ and $g_2^{ij} (u)$ are almost
compatible.
\el

Let us introduce the affinor
\be
v^i_j (u) = g_1^{is} (u) g_{2, sj} (u) \label{aff}
\ee
and consider the Nijenhuis tensor of this affinor
\be
N^k_{ij} (u) = v^s_i (u) {\pa v^k_j \over \pa u^s}
- v^s_j (u) {\pa v^k_i \over \pa u^s}  +
v^k_s (u) {\pa v^s_i \over \pa u^j} -
v^k_s {\pa v^s_j \over \pa u^i},  \label{nij}
\ee
following \cite{[34]}, where were similarly considered
the affinor $v^i_j (u)$
and its Nijenhuis tensor for two flat metrics.

\bt  \label{pochti}
The metrics
$g_1^{ij} (u)$ and $g_2^{ij} (u)$ are almost
compatible if and only if the corresponding Nijenhuis tensor
$N^k_{ij} (u)$ (\ref{nij}) vanishes.
\et

\bl \label{l2}
The following identities are always fulfilled:
\be
g_{1, sp} (u) N^p_{rq} (u) g_2^{ri} (u) g_2^{qj} (u) g_2^{sk} (u)
=M^{kij} (u) + M^{jki} (u) + M^{jik} (u),  \label{mn1}
\ee
\bea
&&
2( M^{jki} (u) + M^{jik} (u))=\nn\\
&&
g_{1, sp} (u) N^p_{rq} (u) g_2^{ri} (u) g_2^{qj} (u) g_2^{sk} (u) +
g_{1, sp} (u) N^p_{rq} (u) g_2^{rk} (u) g_2^{qj} (u) g_2^{si} (u),
\label{mn2}
\eea
\be
2 M^{kij} (u) =
g_{1, sp} (u) N^p_{rq} (u) g_2^{ri} (u) g_2^{qj} (u) g_2^{sk} (u) -
g_{1, sp} (u) N^p_{rq} (u) g_2^{rk} (u) g_2^{qj} (u) g_2^{si} (u).
 \label{mn3}
\ee
\el

\bc
The tensor $M^{ijk} (u)$ vanishes if and only if
the Nijenhuis tensor (\ref{nij}) vanishes.
\ec

In the papers \cite{[35]}--\cite{mokh2}, the present author
studied reductions in the general problem on compatible flat metrics,
connected with the associativity equations, namely, the following
ansatz in formula (\ref{(2.9)}):
$$
h^i (v) = \eta^{is} {\pa \Phi \over \pa v^s},
$$
where
$\Phi (v^1,...,v^N)$ is a function of $N$ variables.

Correspondingly,
in this case the metrics have the form:

\be
g_1^{ij} (v) = \eta^{ij}, \ \ \ g_2^{ij} (v) = \eta^{is}
\eta^{jp} {\pa^2 \Phi \over \pa v^s \pa v^p}.  \label{m1}
\ee

\bt  \label{mo1}
If metrics (\ref{m1}) are almost compatible, then
they are always compatible.
Moreover, in this case, the metric $g_2^{ij} (v)$ is
necessarily also flat, that is, metrics (\ref{m1})
form a flat pencil of metrics. Condition of
almost compatibility of metrics (\ref{m1})
has the form
\be
\eta^{sp} {\pa^2 \Phi \over \pa v^p \pa v^i}
{\pa^3 \Phi \over \pa v^s \pa v^j \pa v^k} =
\eta^{sp} {\pa^2 \Phi \over \pa v^p \pa v^k}
{\pa^3 \Phi \over \pa v^s \pa v^j \pa v^i} \label{m2}
\ee
and coincides with the condition of
compatible deformation of two Frobenius algebras
which was derived and studied by the author in \cite{[36]}--\cite{mokh2}.
\et

In particular, in the author's papers
\cite{[36]}--\cite{mokh2}, it is proved that
in the two-component case ($N=2$), for $\eta^{ij} =
\varepsilon^i \delta^{ij},$ $ \varepsilon^i = \pm 1,$
condition (\ref{m2}) is equivalent to the following
linear second-order partial differential equation with constant
coefficients:
\be
\alpha \left ( \varepsilon^1
{\pa^2 \Phi \over \pa (v^1)^2} -
 \varepsilon^2
{\pa^2 \Phi \over \pa (v^2)^2} \right ) = \beta
 {\pa^2 \Phi \over \pa v^1 \pa v^2},
\ee
where $\alpha$ and $\beta$ are arbitrary constants.

\section{Compatible metrics and Nijenhuis tensor} \label{sect2}

Let us prove the second part of Theorem \ref{tmo1}.
In the previous section it is proved, in particular,
that it always follows from compatibility (and even almost compatibility)
of metrics that the corresponding Nijenhuis tensor vanishes.

Assume that a pair of metrics
$g^{ij}_1 (u)$ and $g^{ij}_2 (u)$ is nonsingular, that is,
the eigenvalues of this pair of the metrics are distinct,
and assume also that the corresponding Nijenhuis tensor vanishes.
Let us prove that, in this case, the metrics
$g^{ij}_1 (u)$ and $g^{ij}_2 (u)$ are compatible
(their almost compatibility follows from the previous
section).

Obviously, that the eigenvalues of the pair of the metrics
$g^{ij}_1 (u)$ and $g^{ij}_2 (u)$  coincide with the
eigenvalues of affinor $v^i_j (u)$. But it is well-known
that if all eigenvalues of an affinor are distinct, then
it always follows from the vanishing of the Nijenhuis tensor
of this affinor that there exist local coordinates such that
in these coordinates
the affinor reduces to a diagonal form in the corresponding neighbourhood
\cite{nij1} (see also \cite{ha}).

So we can consider further that the affinor
$v^i_j (u)$ is diagonal in the local coordinates
 $u^1,...,u^N$,
\be
 v^i_j (u) = \lambda^i (u) \delta^i_j,
\ee
where is no summation over the index $i$,
and $\lambda^i (u), \ i=1,...,N,$ are the eigenvalues
of the pair of metrics $g^{ij}_1 (u)$ and $ g^{ij}_2 (u)$,
which are assumed distinct:
\be
\lambda^i \neq \lambda^j, {\rm \ for\ } i\neq  j.
\ee

\bl              \label{n1}
If the affinor (\ref{aff}) is diagonal in a certain local
coordinates and all its eigenvalues are distinct, then,
in these coordinates, the metrics $g^{ij}_1 (u) $ and
$g^{ij}_2 (u)$ are also necessarily diagonal.
\el

Actually, we have
$$g^{ij}_1 (u) = \lambda^i (u) g^{ij}_2 (u).$$
It follows from symmetry of the metrics $g^{ij}_1 (u)$ and
$ g^{ij}_2 (u)$
that for any $i, j$
\be
(\lambda^i (u) - \lambda^j (u))  g^{ij}_2 (u) = 0,
\ee
where is no summation over indices, that is,
$$ g^{ij}_2 (u) = g^{ij}_1 (u) =0
{\rm \  for \ } i\neq  j.$$

\bl     \label{n2}
Let an affinor $w^i_j (u)$ be diagonal in a certain
local coordinates
$u= (u^1,...,u^N)$, that is, $w^i_j (u) =\mu^i (u) \delta^i_j$.
\bitem
\item[1)] If all the eigenvalues of an diagonal affinor are distinct,
that is,
$\mu^i (u) \neq \mu^j (u)$ for $i \neq  j$, then
the Nijenhuis tensor of this affinor vanishes if and only if
the $i$th eigenvalue $\mu^i (u)$ depend only on the coordinate $u^i.$
\item[2)] If all the eigenvalues coincide, then
the Nijenhuis tensor vanishes.
\item[3)] In the general case of a diagonal affinor,
the Nijenhuis tensor vanishes if and only if
\be
{\pa \mu^i \over \pa u^j} = 0
\ee
for all $i, j$ such that $\mu^i (u) \neq \mu^j (u).$
\eitem
\el

Actually, for any diagonal affinor $w^i_j (u) =\mu^i (u) \delta^i_j,$
the Nijenhuis tensor $N^k_{ij} (u)$ has the form:
$$N^k_{ij} (u) = (\mu^i - \mu^k) {\pa \mu^j \over \pa u^i} \delta^{kj}
- (\mu^j - \mu^k) {\pa \mu^i \over \pa u^j} \delta^{ki}$$
(no summation over indices).
Thus, the Nijenhuis tensor vanishes if and only if
for all $i, j$
$$(\mu^i (u) - \mu^j (u)) {\pa \mu^i \over \pa u^j} =0,$$
where is no summation over indices.

It follows immediately from Lemmas \ref{n1} and \ref{n2}
that for any nonsingular pair of
almost compatible metrics there always exist local coordinates
in which the metrics have the form
$$g^{ij}_2 (u) = g^i (u) \delta^{ij}, \ \ \
g^{ij}_1 (u) = \lambda^i (u^i) g^i (u) \delta^{ij},
\ \ \ \lambda^i = \lambda^i (u^i), \ i=1,...,N.$$

Moreover, we derive immediately that any diagonal
metrics of the form $g^{ij}_2 (u) = g^i (u) \delta^{ij}$ and
$g^{ij}_1 (u) = f^i (u^i) g^i (u) \delta^{ij}$
for any nonzero functions $f^i (u^i),$ $ i=1,...,N,$
(they can be here, for example, coinciding nonzero constants, that is,
the pair of metrics may be ``singular'') are always almost compatible.
We shall prove now that they are always compatible.
Then Theorem \ref{tmo1} will be completely proved.

Consider diagonal metrics
$g^{ij}_2 (u) = g^i (u) \delta^{ij}$ and
$g^{ij}_1 (u) = f^i (u^i) g^i (u) \delta^{ij},$
where $f^i (u^i),$ $i=1,...,N,$ are arbitrary functions of
a single variable, which are not equal to zero identically,
and consider their arbitrary linear combination
$$g^{ij} (u) = (\lambda_2 + \lambda_1 f^i (u^i)) g^i (u) \delta^{ij},$$
where $\lambda_1$ and $\lambda_2$ are arbitrary constants such that
$\det (g^{ij} (u)) \not\equiv 0.$

Let us prove that relation (\ref{kr}) is always fulfilled for
the corresponding tensors of Riemannian curvature.

Recall that for any diagonal metric
$\Gamma^i_{jk} (u) =0$ if all the indices $i, j, k$
are distinct.
Correspondingly,
$R^{ij}_{kl} (u) = 0$
if all the indices $i, j, k, l $  are distinct.
Besides, as a result of the well-known symmetries of
the tensor of Riemannian curvature we have:
$$R^{ii}_{kl} (u) = R^{ij}_{kk} (u) =0,$$
$$R^{ij}_{il} (u) = -R^{ij}_{li} (u) = R^{ji}_{li} (u) =
- R^{ji}_{il} (u).$$

Thus, it is sufficient to prove relation (\ref{kr}) only
for the following components of the tensor of Riemannian curvature:
 $R^{ij}_{il} (u)$,
where $i \neq j,$ $\ i \neq  l$.

For any diagonal metric
$g^{ij}_2 (u) = g^i (u) \delta^{ij}$   we have
$$\Gamma^i_{2, ik} (u) = \Gamma^i_{2, ki} (u) =
- {1 \over 2 g^i (u)} {\pa g^i \over \pa u^k}, \ \ \
{\rm \ for\  any\ } i, k; $$
$$\Gamma^i_{2, jj} (u) = {1 \over 2} {g^i (u) \over
(g^j (u))^2 } {\pa g^j \over \pa u^i},\ \ \ i \neq j.$$

\bea
&&
R^{ij}_{2, il} (u) = g^i (u) R^j_{2, iil} (u) = \nn\\
&&
g^i (u) \left ( - {\pa \Gamma^j_{2, il} \over \pa u^i}  +
{\pa \Gamma^j_{2, ii}  \over \pa u^l} - \sum_{s=1}^N
\Gamma^j_{2, si} (u) \Gamma^s_{2, il} (u) +
\sum_{s=1}^N \Gamma^j_{2, sl} (u) \Gamma^s_{2, ii} (u) \right ).
\eea

It is necessary to consider separately two different cases.

1) $j \neq l$.

\bea
&&
R^{ij}_{2, il} (u) =
g^i (u) \left (
{\pa \Gamma^j_{2, ii}  \over \pa u^l} -
\Gamma^j_{2, ii} (u) \Gamma^i_{2, il} (u) +
\Gamma^j_{2, jl} (u) \Gamma^j_{2, ii} (u)
 +
\Gamma^j_{2, ll} (u) \Gamma^l_{2, ii} (u) \right )
=  \nn\\
&&
{1 \over 2} g^i (u)  {\pa \over \pa u^l}
\left ( {g^j (u) \over (g^i (u))^2 } {\pa g^i \over \pa u^j }\right )
 + {1 \over 4} {g^j (u) \over (g^i (u))^2 }
{\pa g^i \over \pa u^j} {\pa g^i \over \pa u^l} \nn\\
&&
- {1 \over 4 g^i (u)}
{\pa g^i \over \pa u^j} {\pa g^j \over \pa u^l}
+ {1 \over 4} {g^j (u) \over g^i (u) g^l (u) }
{\pa g^l \over \pa u^j} {\pa g^i \over \pa u^l}.\label{r1}
\eea

Respectively, for the metric
$$g^{ij} (u) = (\lambda_2 + \lambda_1 f^i (u^i)) g^i (u) \delta^{ij},$$
we obtain (we use here that all the indices $i, j, l$ are distinct):
\bea
&&
R^{ij}_{il} (u)=
(\lambda_2 + \lambda_1 f^j (u^j)) \left [
{1 \over 2} g^i (u)  {\pa \over \pa u^l}
\left ( {g^j (u) \over (g^i (u))^2 } {\pa g^i \over \pa u^j }\right )
 + {1 \over 4} {g^j (u) \over (g^i (u))^2 }
{\pa g^i \over \pa u^j} {\pa g^i \over \pa u^l} \right.\nn\\
&&
- \left. {1 \over 4 g^i (u)}
{\pa g^i \over \pa u^j} {\pa g^j \over \pa u^l}
+ {1 \over 4} {g^j (u) \over g^i (u) g^l (u) }
{\pa g^l \over \pa u^j} {\pa g^i \over \pa u^l} \right ]=
\lambda_1 R^{ij}_{1, il} (u) + \lambda_2 R^{ij}_{2, il} (u). \label{rr1}
\eea

2) $j=l$.

\bea
&&
R^{ij}_{2, ij} (u) =
g^i (u) \left (
- {\pa \Gamma^j_{2, ij}  \over \pa u^i} +
{\pa \Gamma^j_{2, ii}  \over \pa u^j} -
\Gamma^j_{2, ii} (u) \Gamma^i_{2, ij} (u) - \right. \nn\\
&&
\left. \Gamma^j_{2, ji} (u) \Gamma^j_{2, ij} (u)
 +  \sum_{s=1}^N
\Gamma^j_{2, sj} (u) \Gamma^s_{2, ii} (u) \right )
=  \nn\\
&&
{1 \over 2} g^i (u)  {\pa \over \pa u^i}
\left ( {1 \over g^j (u) } {\pa g^j \over \pa u^i }\right )
+ {1 \over 2} g^i (u)  {\pa \over \pa u^j}
\left ( {g^j (u) \over (g^i (u))^2 } {\pa g^i \over \pa u^j }\right )
+ {1 \over 4} {g^j (u) \over (g^i (u))^2 }
{\pa g^i \over \pa u^j} {\pa g^i \over \pa u^j} \nn\\
&&
- {1 \over 4} {g^i (u) \over (g^j(u))^2}
{\pa g^j \over \pa u^i} {\pa g^j \over \pa u^i}
+ {1 \over 4 g^j (u)}
{\pa g^j \over \pa u^i} {\pa g^i \over \pa u^i} -
\sum_{s \neq  i} {1 \over 4} {g^s (u) \over g^i (u) g^j (u) }
{\pa g^j \over \pa u^s} {\pa g^i \over \pa u^s}.\label{r2}
\eea

Respectively, for the metric
$$g^{ij} (u) = (\lambda_2 + \lambda_1 f^i (u^i)) g^i (u) \delta^{ij},$$
we obtain (we use here that the indices $i, j$ are distinct):

\bea
&&
R^{ij}_{ij} (u) =
{1 \over 2} (\lambda_2 + \lambda_1 f^i (u^i))
g^i (u)  {\pa \over \pa u^i}
\left ( {1 \over g^j (u) } {\pa g^j \over \pa u^i }\right )
+ \nn\\
&&
{1 \over 2} g^i (u)  {\pa \over \pa u^j}
\left ( {
(\lambda_2 + \lambda_1 f^j (u^j))
g^j (u) \over (g^i (u))^2 } {\pa g^i \over \pa u^j }\right )
+ {1 \over 4}
(\lambda_2 + \lambda_1 f^j (u^j))
{g^j (u) \over (g^i (u))^2 }
{\pa g^i \over \pa u^j} {\pa g^i \over \pa u^j} \nn\\
&&
- {1 \over 4}
(\lambda_2 + \lambda_1 f^i (u^i))
 {g^i (u) \over (g^j(u))^2}
{\pa g^j \over \pa u^i} {\pa g^j \over \pa u^i}
+ {1 \over 4 g^j (u)}
{\pa g^j \over \pa u^i} {\pa
((\lambda_2 + \lambda_1 f^i (u^i))
g^i) \over \pa u^i} -
\nn\\
&&
{1 \over 4 g^i (u)}
{\pa g^i \over \pa u^j} {\pa
((\lambda_2 + \lambda_1 f^j (u^j)))
g^j \over \pa u^j} -
\sum_{s \neq  i, \ s \neq  j}
{1 \over 4} {
(\lambda_2 + \lambda_1 f^s (u^s))
g^s (u) \over g^i (u) g^j (u) }
{\pa g^j \over \pa u^s} {\pa g^i \over \pa u^s}
=\nn\\
&&
\lambda_1 R^{ij}_{1, ij} (u) + \lambda_2 R^{ij}_{2, ij} (u).
\eea

Theorem \ref{tmo1} is proved.
Moreover, the complete explicit description
of nonsingular pairs of compatible and almost compatible
metrics is obtained and the following theorem is proved:

\bt   \label{teom}
An arbitrary nonsingular pair of metrics is compatible
if and only if there exist local coordinates $u = (u^1,...,u^N)$
such that
$g^{ij}_2 (u) = g^i (u) \delta^{ij}$ and
$g^{ij}_1 (u) = f^i (u^i) g^i (u) \delta^{ij},$
where $f^i (u^i),$ $i=1,...,N,$ are arbitrary functions
of single variables (of course, in the case of nonsingular pair of metrics,
these functions are not equal to each other if they are constants and they
are not equal identically to zero).
\et

Let us consider here the problem on nonsingular pairs
of compatible flat metrics. It follows from Theorem \ref{teom}
that it is
sufficient to classify flat metrics of the form
$g_2^{ij} (u) = g^i (u) \delta^{ij}$ and
$g_1^{ij} (u) = f^i (u^i) g^i (u) \delta^{ij},$
where $f^i (u^i),$ $i= 1,...,N,$ are arbitrary functions of
single variables.

The  problem of description of diagonal flat metrics,
that is, flat metrics
$g_2^{ij} (u) = g^i (u) \delta^{ij},$
is a classical problem of differential geometry.
This problem is equivalent to the problem
of description of curvilinear orthogonal coordinate systems in
a pseudo-Euclidean space and it was studied in detail and mainly solved
in the beginning of the 20th century (see \cite{da}).
Locally, such coordinate systems are determined by
$n(n-1)/2$ arbitrary functions of two variables (see \cite{ca},
\cite{bi}). Recently, Zakharov showed that the Lam\'{e}
equations describing curvilinear orthogonal coordinate systems
can be integrated by the inverse scattering method \cite{za}
(see also an algebraic-geometric approach in \cite{kri}).
The condition that the metric
$g_1^{ij} (u) = f^i (u^i) g^i (u) \delta^{ij}$ is also flat
gives exactly $n(n-1)/2$ equations which are linear with respect to
the functions $f^i (u^i)$.  Note that, in this case,
components (\ref{r1})
of tensor of Riemannian curvature
automatically vanish as a result of formula (\ref{rr1}).
And the vanishing of components (\ref{r2}) gives
the corresponding $n(n-1)/2$ equations. In particular, in the case
$N=2$ this completely solves the problem of description
of nonsingular pairs of compatible flat metrics.
In the next section we give their complete description.
But without a doubt this problem is integrable in general for
any $N$. We are going to devote to this question another work.
In particular, it is very interesting to classify
all the $n$-orthogonal curvilinear coordinate systems in a
pseudo-Euclidean space, for which the functions
$f^i (u^i) = (u^i)^n$ define compatible flat metrics
(respectivly, separately for $n=1$; $n=1, 2$; $n=1, 2, 3,$ and so on).

\section{Two-component compatible flat metrics} \label{sect5}

We present here the complete description of
nonsingular pairs of two-component compatible
flat metrics (see also \cite{[35]}, \cite{mokh1},
\cite{mokh2}, where an integrable
 homogeneous system of hydrodynamic type,
describing all the two-component compatible flat metrics, was derived and
investigated).

It is proved above that for any nonsingular pair
of two-component compatible metrics $g^{ij}_1 (u)$
and $g^{ij}_2 (u)$ there always exist local coordinates
$u^1,...,u^N$ such that
\be
(g^{ij}_2 (u) ) = \left ( \begin{array}{cc}
{\varepsilon^1 \over (b^1 (u))^2 } & 0\\
0 & {\varepsilon^2 \over (b^2 (u))^2 }
\end{array} \right ), \ \ \ \ \
(g^{ij}_1 (u) ) = \left ( \begin{array}{cc}
{\varepsilon^1  f^1 (u^1) \over (b^1 (u))^2 } & 0\\
0 & {\varepsilon^2  f^2 (u^2) \over (b^2 (u))^2 }
\end{array} \right ), \label{metr}
\ee
where $\varepsilon^i = \pm 1, \ i=1, 2;$
$b^i (u)$ and $f^i (u^i), \ i=1, 2,$ are arbitrary nonzero functions
of the corresponding single variables.

\bl \label{moflat}
An arbitrary diagonal metric $g^{ij}_2 (u)$ (\ref{metr})
is flat if and only if the functions $b^i (u),$  $ i=1, 2,$
are solutions of the following linear system:
\be
{\pa b^2 \over \pa u^1} = \varepsilon^1
{\pa F \over \pa u^2} b^1 (u),\ \ \ \
{\pa b^1 \over \pa u^2} =
- \varepsilon^2 {\pa F \over \pa u^1} b^2 (u),
\label{sys}
\ee
where $F(u)$ is an arbitrary function.
\el

\bt
The metrics $g^{ij}_1 (u)$ and $g^{ij}_2 (u)$ (\ref{metr})
form a flat pencil of metrics if and only if
the functions $b^i (u), \ i=1, 2,$ are solutions of
the linear system (\ref{sys}), where the function $F(u)$ is
a solution of the following linear equation:
\be
2 {\pa^2 F \over \pa u^1 \pa u^2} (f^1 (u^1) -  f^2 (u^2))
+ {\pa F \over \pa u^2} {d f^1 (u^1)  \over d u^1} -
{\pa F \over \pa u^1} {d f^2 (u^2) \over d u^2} =0. \label{lequa}
\ee
\et

In the case, if the eigenvalues of the pair of the metrics
$g^{ij}_1 (u)$ and $g^{ij}_2 (u)$ are not only distinct but
also are not constants, we can always choose local coordinates
such that $f^1 (u^1) = u^1,$  $f^2 (u^2) = u^2$ (see also remark in
\cite{[34]}).
In this case, equation (\ref{lequa}) has the form
\be
2 {\pa^2 F \over \pa u^1 \pa u^2} (u^1 -  u^2)
+ {\pa F \over \pa u^2} -
{\pa F \over \pa u^1} = 0. \label{lequa2}
\ee

Let us continue this recurrent procedure for the metrics
$G^{ij}_{n+1} (u) = v^i_j (u) G^{ij}_n(u)$ with the help of
the affinor $v^i_j (u) = u^i \delta^i_j.$

\bt
Three metrics
\be
(G^{ij}_n (u) ) = \left ( \begin{array}{cc}
{\varepsilon^1  (u^1)^n \over (b^1 (u))^2 } & 0\\
0 & {\varepsilon^2   (u^2)^n \over (b^2 (u))^2 }
\end{array} \right ), \ \ \ \ \ n=0, 1, 2,\label{metr2}
\ee
form a flat pencil of metrics (pairwise compatible)
if and only if the functions
$b^i (u), \ i=1, 2,$ are solutions of the linear system (\ref{sys}),
where
\be
F(u) = c \ln (u^1 - u^2),
\ee
$c$ is an arbitrary constant.
Already the metric $G^{ij}_3 (u)$ is flat only in the
most trivial case, when $c =0,$
and, respectively, $b^1 = b^1 (u^1),$ $b^2 = b^2 (u^2)$.

The metric $G^{ij}_3 (u)$ is a metric of
nonzero constant Riemannian curvature $K \neq 0$
if and only if
\be
(b^1 (u))^2 = (b^2 (u))^2 = {\varepsilon^2 \over 4 K}
(u^1 - u^2),\ \ \ \ \varepsilon^1 = - \varepsilon^2,
\ \ \ \ c = \pm {1 \over 2}.
\ee
\et

\section {Almost compatible metrics which are not compatible}

\bl   \label{molem}
Two-component diagonal conformally Euclidean metric
$$g^{ij} (u) =  \exp (a (u)) \delta^{ij}, \
1 \leq i,j \leq 2,$$
is flat if and only if
the function $a(u)$ is harmonic, that is,
\be
\Delta a \equiv {\pa^2 a \over \pa (u^1)^2} +
{\pa^2 a \over \pa (u^2)^2} =0.
\ee
\el

In particular, the metric
$g_1^{ij} (u) =  \exp (u^1 u^2) \delta^{ij}, \
1 \leq i,j \leq 2$, is flat. Obviously, that the flat metrics
$g_1^{ij} (u) =  \exp (u^1 u^2) \delta^{ij}, \
1 \leq i,j \leq 2$, and $g_2^{ij} (u) = \delta^{ij}, \ 1 \leq i, j \leq 2,$
are almost compatible, for them the Nijenhuis tensor (\ref{nij}) vanishes.
But it follows from Lemma \ref{molem} that these metrics are not
compatible, their sum is not a flat metric.

Similarly it is possible to construct also other counterexamples to
Theorem \ref{tfer}. Moreover, the following statement is true.

\bp
Any nonconstant real harmonic function $a(u)$
defines a pair of almost compatible metrics
 $g_1^{ij} (u) =  \exp (a (u)) \delta^{ij}, \
1 \leq i,j \leq 2$, and $g_2^{ij} (u) = \delta^{ij}, \ 1 \leq i, j \leq 2,$
which are not compatible. These metrics are compatible if
and only if $a = a(u^1 \pm i u^2).$
\ep

Let us construct also almost compatible metrics of
constant Riemannian curvature, which are not compatible.

\bl   \label{molem2}
Two-component diagonal conformally Euclidean metric
$$g^{ij} (u) =  \exp (a (u)) \delta^{ij}, \
1 \leq i,j \leq 2,$$
is a metric of constant Riemannian curvature $K$
if and only if the function $a(u)$ is a solution of the Liouville equation
\be
\Delta a \equiv {\pa^2 a \over \pa (u^1)^2} +
{\pa^2 a \over \pa (u^2)^2} = 2 K e^{- a (u)}.  \label{liuv}
\ee
\el

\bp
For the metrics
$g_1^{ij} (u) =  \exp (a (u)) \delta^{ij}, \
1 \leq i,j \leq 2$, and $g_2^{ij} (u) = \delta^{ij}, \ 1 \leq i, j \leq 2,$
the corresponding Nijenhuis tensor vanishes, that is, they are
always almost compatible.
But they are compatible metrics of constant Riemannian curvature
$K \neq 0$ and $0$, respectively, if and only if the function
$a(u)$ is constant.
\ep

Note, that all the one-component ``metrics'' are always
compatible, and all the one-component local Poisson structures
of hydrodynamic type are also always compatible.
Let us construct for any $N > 1$ examples of
almost compatible metrics which are not compatible.

\bp
The metrics
$g_1^{ij} (u) =  b (u) \delta^{ij}, \
1 \leq i,j \leq N$ and
$g_2^{ij} (u) = \delta^{ij}, \ 1 \leq i, j \leq N,$
where $b (u)$ is an arbitrary function, are always almost compatible,
the corresponding Nijenhuis tensor vanishes.
But they are compatible real metrics only in the most trivial
case, when the function $b(u)$ is constant.
Complex metrics are compatible if and only if either the function
 $b(u)$ is constant, or  $N=2$ and $b(u) = b (u^1 \pm i u^2)$.
\ep

\section{Compatible metrics and nonlocal Poisson structures
of hydrodynamic type}

Nonlocal Poisson structures of hydrodynamic type were introduced
and studied in the work of the present author and Ferapontov
\cite{[38]}
(see also \cite{[40]}--\cite{[39]}). They have the following form:

\be
\{ I,J \} = \int {\delta I \over \delta u^i(x)}
\left ( g^{ij}(u(x)) {d \over dx} + b^{ij}_k (u(x))\, u^k_x
+ K u^i_x \left ( {d \over dx} \right )^{-1} u^j_x \right )
{\delta J \over \delta u^j(x)} dx,
\label{nonl}
\ee
where $K$ is an arbitrary constant.

The bracket of form (\ref{nonl}) is called {\it nondegenerate} if
$\det (g^{ij}) (u) \not\equiv 0$.

\bt [\cite{[38]}] \label{mofer}
If $\det (g^{ij}) (u) \not\equiv 0$, then bracket (\ref{nonl})
is a Poisson bracket, that is, it is skew-symmetric and
satisfies the Jacobi identity, if and only if
\bitem
\item [(1)] $g^{ij} (u)$ is an arbitrary pseudo-Riemannian
contravariant metric of constant Riemannian curvature $K$,

\item [(2)] $b^{ij}_k (u) = - g^{is} (u) \Gamma ^j_{sk} (u),$ where
$\Gamma^j_{sk} (u)$ is the Riemannian connection generated
by the contravariant metric $g^{ij} (u)$
 (the Levi--Civita connection).
\eitem
\et

\bp
Nonlocal nondegenerate Poisson brackets of form (\ref{nonl})
are compatible if and only if their metrics are compatible.
\ep

In \cite{[40]} Ferapontov introduced and studied more general
nonlocal Poisson brackets of hydrodynamic type, namely, the brackets
of the following form:
\bea
&& \
\{ I,J \} = \int   {\delta I \over \delta u^i(x) }
\left ( g^{ij}(u(x)) {d \over dx} + b^{ij}_k (u(x))\, u^k_x \right.
\nn\\
&& \
\left. + \sum_{\alpha =1}^L
(w^{\alpha})^i_k (u)
u^k_x \left ( {d \over dx} \right )^{-1}
(w^{\alpha})^j_s (u) u^s_x \right )
{\delta J \over \delta u^j(x)} dx,
\ \ \ \ \det (g^{ij}) (u) \not\equiv 0.
\label{nonl2}
\eea

\bt [\cite{[40]}] \label{ferap}
Bracket (\ref{nonl2}) is a Poisson bracket, that is,
it is skew-symmetric and satisfies the Jacobi identity,
if and only if
\bitem
\item [(1)] $b^{ij}_k (u) = - g^{is} (u) \Gamma ^j_{sk} (u),$ where
$\Gamma^j_{sk} (u)$ is the Riemannian connection generated by the
contravariant metric $g^{ij} (u)$
 (the Levi--Civita connection),

\item [(2)] the pseudo-Riemannian
metric $g^{ij} (u)$ and the set of affinors
$(w^{\alpha})^i_j (u)$ satisfy the relations:

\be
g_{ik} (u) (w^{\alpha})^k_j (u) =
g_{jk} (u) (w^{\alpha})^k_i (u),
\ \ \ \alpha = 1,...,L,  \label{peter1}
\ee
\be
\nabla_k (w^{\alpha})^i_j (u) =
\nabla_j (w^{\alpha})^i_k (u),
\ \ \ \alpha = 1,...,L,   \label{peter2}
\ee
\be
R^{ij}_{kl} (u) =  \sum_{\alpha =1}^L
\left ( (w^{\alpha})^i_k (u) (w^{\alpha})^j_l (u)
- (w^{\alpha})^j_k (u) (w^{\alpha})^i_l (u) \right ). \label{gauss}
\ee
Moreover, the family of affinors $w^{\alpha} (u)$ is commutative:
$[w^{\alpha}, w^{\beta}] =0.$
\eitem
\et

\bp
If nonlocal Poisson brackets of form (\ref{nonl2})
are compatible, then their metrics are also compatible.
\ep

Actually,
if nonlocal Poisson brackets of form (\ref{nonl2})
are compatible, then it follows from the conditions of
compatibility and from Theorem \ref{ferap} that, first,
relation (\ref{sv}) is fulfilled and, secondly, the curvature tensor
for the metric
$g^{ij} (u) = \lambda_1 g^{ij}_1 (u) + \lambda_2 g^{ij}_2 (u)$
has the form
\bea
&&
R^{ij}_{kl} (u) =
\sum_{\alpha =1}^{L_1} \lambda_1
\left ( (w_1^{\alpha})^i_k (u) (w_1^{\alpha})^j_l (u)
- (w_1^{\alpha})^j_k (u) (w_1^{\alpha})^i_l (u) \right ) + \nn\\
&&
\sum_{\alpha =1}^{L_2} \lambda_2
\left ( (w_2^{\alpha})^i_k (u) (w_2^{\alpha})^j_l (u)
- (w_2^{\alpha})^j_k (u) (w_2^{\alpha})^i_l (u) \right ) =
\lambda_1 R^{ij}_{1, kl} (u) + \lambda_2 R^{ij}_{2, kl} (u).\nn
\eea

Apparently, the converse statement is also always true
(this is only our conjecture which is not strictly proved in the
most general case at this moment).

Relations (\ref{peter1})--(\ref{gauss}) are nothing but
the Gauss--Peterson--Codazzi equations for $N$-di\-men\-si\-o\-nal
surfaces $M$ with flat normal connections in a pseudo-Euclidean
space $E^{N+L}$. Here $g_{ij} (u)$ is
the first fundamental form of the surface $M$, and $w^{\alpha} (u)$
are the Weingarten operators \cite{[40]}. We consider more in detail
the case of compatible nonlocal Poisson structures of
hydrodynamic type that correspond to surfaces with holonomic
net of curvature lines (see \cite{[34]}).
This case is the most interesting for applications
(here we do not give numerous important examples, see, for example,
in \cite{[40]}, \cite{[28]},
\cite{al} or in the author's survey \cite{mokh4}).

\bp
Let two nonlocal Poisson brackets of form (\ref{nonl2}) correspond to
surfaces with holonomic net of curvature lines and be given
in coordinates of curvature lines. In this case,
if the corresponding pair of metrics is nonsingular, then
the nonlocal Poisson structures are compatible if and only if
their metrics are compatible.
\ep

In this case the metrics $g^{ij}_1 (u) = g^i_1 (u) \delta^{ij}$
and $g^{ij}_2 (u) = g^i_2 (u) \delta^{ij}$,
and also the Weingarten operators $(w_1^{\alpha})^i_j (u)
= (w^{\alpha}_1)^i (u) \delta^i_j$
and $w_2^{\alpha} (u)
= (w^{\alpha}_2)^i (u) \delta^i_j$  are diagonal in
the coordinates under consideration.
For any such ``diagonal'' case,
condition (\ref{peter1}) is automatically fulfilled,
all the Weingarten operators commute, conditions (\ref{peter2})
and (\ref{gauss}) have the following form, respectively:
\be
2 g^i (u) {\pa (w^{\alpha})^i \over \pa u^k} =
((w^{\alpha})^i - (w^{\alpha})^k) {\pa g^i \over \pa u^k}
{\rm \  for \ all \ } i \neq k,\label{hol1}
\ee
\be
R^{ij}_{ij} (u) = \sum_{\alpha =1}^L (w^{\alpha})^i (u)
(w^{\alpha})^j (u). \label{hol2}
\ee

It follows from nonsingularity of the pair of the metrics
and from compatibility of the metrics that the corresponding
Nijenhuis tensor vanishes and there exist functions
 $f^i (u^i), i=1,...,N,$ such that:
$$g^i_1 (u) = f^i (u^i) g^i_2 (u).$$
Using relations (\ref{hol1}) and (\ref{hol2}), it is easy to prove that
in this case it follows from compatibility of the metrics that
an arbitrary linear combination of nonlocal Poisson
brackets under consideration is also a Poisson bracket.


\medskip

\begin{flushleft}
Centre for Nonlinear Studies,\\
L.D.Landau Institute for Theoretical Physics,\\
Russian Academy of Sciences,\\
ul. Kosygina, 2,\\
Moscow, 117940  Russia\\
e-mail: mokhov@genesis.mi.ras.ru, mokhov@landau.ac.ru\\

\smallskip

Department of Mathematics,\\
University of Paderborn,\\
Paderborn, Germany\\
e-mail: mokhov@uni-paderborn.de
\end{flushleft}

\end{document}